arXiv:math/9703202v1 [math.GR] 25 Mar 1997

# DUALITY AND LOCAL GROUP COHOMOLOGY

P R Hewitt

Univ of Toledo

23 March 97

Abstract. Recently, Meierfrankenfeld has published three theorems on the cohomology of a finitary module. They cover the local determination of complete reducibility; the local splitting of group extensions; and the representation of locally split extensions in the double dual. In this note we derive all three by combining a certain duality between homology and cohomology with the continuity of homology.

## Local Cohomology

We describe a general framework for three recent results of Meierfrankenfeld's [M], on the cohomology of finitary modules. It turns out that Theorem A below — which began life as an attempt to understand the third of this trio — includes as corollaries all three. Moreover, it implies some new corollaries, and covers all cohomological degrees.

Throughout, $G$ is a group, $k$ is a commutative field, and $\mathcal{L}$ is a *local system* — an upwardly directed collection of subgroups whose union is $G$. Modules are $k$-spaces. If $V$ is a module, then $V^\vee$ denotes its dual, and $\langle \mid \rangle$ denotes their pairing.

We will explore the relationship between the cohomology of $G$, $\mathrm{H}^*(G, V)$, with the local cohomology, $\mathrm{H}^*(\mathcal{L}, V)$. This latter is defined by taking the limit, with respect to restriction maps, of the cohomology groups $\mathrm{H}^*(L, V)$, for $L$ in $\mathcal{L}$:

$$\mathrm{H}^*(\mathcal{L}, V) := \varprojlim \mathrm{H}^*(L, V).$$

We can define analogous limits, $\mathrm{C}^*(\mathcal{L}, V)$, $\mathrm{Z}^*(\mathcal{L}, V)$, and $\mathrm{B}^*(\mathcal{L}, V)$, for cochains, cocycles, and coboundaries, respectively. What we find is that there are canonical isomorphisms $\mathrm{C}^*(\mathcal{L}, V) \cong \mathrm{C}^*(G, V)$ and $\mathrm{Z}^*(\mathcal{L}, V) \cong \mathrm{Z}^*(G, V)$. Indeed, the localization $\mathrm{C}^*(G, V) \to \mathrm{C}^*(\mathcal{L}, V)$ is defined by restriction to the local subgroups, $\phi \mapsto \{\phi|_L\}$. In the reverse direction, the map $\mathrm{C}^*(\mathcal{L}, V) \to \mathrm{C}^*(G, V)$ splices an $\mathcal{L}$-sequence of $n$-cochains $\{\phi_L\}$ into a cochain on $G$: $g_1, \dots, g_n \mapsto \phi_L(g_1, \dots, g_n)$, where $L$ is any member of $\mathcal{L}$ that contains all the $g_i$. This is well-defined, since in the inverse limit $\mathrm{C}^*(\mathcal{L}, V)$, $\phi_L|_{L'} = \phi_{L'}$, whenever $L \supset L'$.

These maps respect coboundary, and so they induce isomorphisms for the cocycle groups. However, splicing an $\mathcal{L}$-sequence of coboundaries need not yield a coboundary in $G$. In general we obtain only an embedding $\mathrm{B}^*(G, V) \hookrightarrow \mathrm{B}^*(\mathcal{L}, V)$.

Now consider the localization $\mathrm{H}^*(G,V) \to \mathrm{H}^*(\mathcal{L},V)$. An element of $\mathrm{H}^i(\mathcal{L},V)$ is an inverse system of affine flats $\{\phi_L + \mathrm{B}^i(L,V)\}$, where $\phi_L|_{L'} \equiv \phi_{L'} \bmod \mathrm{B}^i(L',V)$ whenever $L \supset L'$. The image of $\mathrm{H}^i(G,V)$ consists of those whose inverse limit is nonempty. Similarly, an element in the kernel of the map $\mathrm{H}^{i+1}(G,V) \to \mathrm{H}^{i+1}(\mathcal{L},V)$ is an inverse system of affine flats $\{\psi_L + \mathrm{Z}^i(L,V)\}$, where $\delta\psi_L|_{L'} = \delta\psi_{L'}$ whenever $L \supset L'$ — viewed modulo the collection of those inverse systems that have a nonempty limit.

The simple 'compactness' proof of [M, Thm 2] applies in arbitrary cohomological degree, to give the following generalization.

**Proposition.** *Let $i$ be a positive integer.*

*a. If $\dim_k \mathrm{B}^i(L,V) < \infty$ for all $L$ in $\mathcal{L}$, then localization $\mathrm{H}^i(G,V) \to \mathrm{H}^i(\mathcal{L},V)$ is surjective in degree $i$.*

*b. If moreover $\dim_k \mathrm{Z}^i(L,V) < \infty$ for all $L$ in $\mathcal{L}$, then additionally localization $\mathrm{H}^{i+1}(G,V) \to \mathrm{H}^{i+1}(\mathcal{L},V)$ is injective in degree $i+1$.* □

We omit the proof because we will not be using this result. The problem with this generalization is that it concerns finitary modules only for $\mathrm{H}^2(G,V)$.

We will write $\mathrm{Ext}^*_G(X,Y)$ instead of $\mathrm{H}^*(G,\mathrm{Hom}_k(X,Y))$, and $\mathrm{Tor}^G_*(X,Y)$ instead of $\mathrm{H}_*(G,X\otimes_k Y)$. In addition to local cohomology, we could also define local homology: $\mathrm{Tor}^{\mathcal{L}}_*(X,Y) = \varinjlim \mathrm{Tor}^L_*(X,Y)$. However, this yields nothing new.

**Lemma 1.** *If $X$ and $Y$ are $G$-modules, then $\mathrm{Tor}^{\mathcal{L}}_*(X,Y) = \mathrm{Tor}^G_*(X,Y)$.*

*Proof.* First of all, note that if $P_* = \mathrm{C}_*(G,\mathbb{Z})$, and $V$ is any $G$-module, then $\varinjlim P_* \otimes_L V = P_* \otimes_G V$. To see this, observe that $P_* \otimes_G V = (P_* \otimes V)/[P_* \otimes V, G]$ and $[P_* \otimes V, G] = \sum_L [P_* \otimes V, L]$.

Since the maps $P_* \otimes_L V \to P_* \otimes_G V$ are surjective, we obtain that the boundaries of the limit are limits of boundaries. Now if an element maps to 0 in a direct limit then it maps to 0 at some finite stage. This tells us that the cycles of the limit are limits of cycles. □

Although local cohomology can differ from global cohomology, we will see that they are the same for the dual of a module.

## Main Theorem

In this section we first record a duality between homology and cohomology. We then exploit this duality, together with the continuity of homology, to obtain a theorem on the local detection of cohomology, valid in all degrees.

**Lemma 2.** *If $f\colon V \to V$ is a linear endomorphism, then*

$$\begin{aligned} \mathrm{ann}_V(\mathrm{im}(f^\vee)) &= \ker(f) \qquad & \mathrm{ann}_{V^\vee}(\mathrm{im}(f)) &= \ker(f^\vee) \\ \mathrm{ann}_V(\ker(f^\vee)) &= \mathrm{im}(f) \qquad & \mathrm{ann}_{V^\vee}(\ker(f)) &= \mathrm{im}(f^\vee) \end{aligned}$$

*Proof.* The first two equalities follow from the identity $\langle f(x) \mid \lambda\rangle = \langle x \mid f^\vee(\lambda)\rangle$. The third follows from the first, since any subspace of $V$ is the annihilator of its annihilator. The last follows from composition of the canonical isomorphisms $\mathrm{ann}_{V^\vee}(\ker(f)) = \mathrm{im}(f)^\vee = \mathrm{im}(f^\vee)$. □

The following is a special case of [B, Prop 2.8.5], but we include a direct proof.

**Lemma 3.** *If $X$ and $Y$ are $L$-modules, then* $\mathrm{Ext}^*_L(X, Y^\vee) = \mathrm{Tor}^L_*(X, Y)^\vee$.

*Proof.* Let $P_* = \mathrm{C}_*(G, \mathbb{Z})$, and note the canonical isomorphism

$$\mathrm{Hom}(P_*, \mathrm{Hom}_k(X, Y^\vee)) = (P_* \otimes X \otimes_k Y)^\vee.$$

This isomorphism can be described as follows: a functional $\lambda$ in $(P_* \otimes X \otimes_k Y)^\vee$ is paired to an additive map $f\colon P_* \to \mathrm{Hom}_k(X, Y^\vee)$ when they satisfy the relation

$$\langle y \mid f(c)(x)\rangle = \langle c \otimes x \otimes y \mid \lambda\rangle.$$

Now if we take fixed points for $L$ in this isomorphism, we obtain that

$$\mathrm{Hom}_L(P_*, \mathrm{Hom}_k(X, Y^\vee)) = (P_* \otimes_L (X \otimes_k Y))^\vee.$$

To finish, apply Lemma 2 to the boundary of $P_* \otimes_L (X \otimes_k Y)$. □

**Theorem A.** *If $X$ and $Y$ are $G$-modules, then localization*

$$\mathrm{Ext}^*_G(X, Y^\vee) \to \mathrm{Ext}^*_{\mathcal{L}}(X, Y^\vee)$$

*is an isomorphism.*

*Proof.* Apply Lemmas 1 and 3, and the fact that the dual of a direct limit of $k$-spaces is the inverse limit of the duals of those spaces. □

**Corollary 1.** *Locally trivial coclasses are trivial in the double dual.*

*Proof.* Extension of scalars from $Y$ to $Y^{\vee\vee}$ factors through localization:

$$\begin{array}{ccc} \mathrm{Ext}^*_G(X, Y) & \to & \mathrm{Ext}^*_G(X, Y^{\vee\vee}) \\ \downarrow & & \| \\ \mathrm{Ext}^*_{\mathcal{L}}(X, Y) & \to & \mathrm{Ext}^*_{\mathcal{L}}(X, Y^{\vee\vee}). \end{array}$$ □

**Corollary 2.** *If $X$ is an arbitrary $G$-module and $Y$ is a finite-dimensional $G$-module, then* $\mathrm{Ext}^*_G(X, Y)$ *is determined locally.* □

In the next section we derive several more corollaries of this theorem, including all three theorems from [M].

## Finitary Cohomology

We now turn to the cohomology of finitary modules. The following characterization of finitary groups — due to Meierfrankenfeld — is the key to understanding their cohomology.

**Lemma 4.** *The finitary group on $V$ is exactly the centralizer of $V^{\vee\vee}/V$.*

*Proof.* Let $x$ be a transformation of $V$. Apply Lemma 2 to $x - 1$ twice: between $V$ and $V^\vee$, and also between $V^\vee$ and $V^{\vee\vee}$. Next, use the fact that a space equals its double dual if and only if it is finite dimensional. □

In light of this characterization, consider the long exact sequence associated to the extension $V \hookrightarrow V^{\vee\vee}$:

$$(\dagger) \qquad \cdots \to \mathrm{Ext}^{i-1}_G(U, V^{\vee\vee}/V) \to \mathrm{Ext}^i_G(U, V) \to \mathrm{Ext}^i_G(U, V^{\vee\vee}) \to \cdots.$$

**Theorem B.** *If $U$ is an arbitrary $G$-module and $V$ is a finitary $G$-module, then the locally trivial portion of* $\operatorname{Ext}^i_G(U,V)$ *lies in the image of* $\operatorname{Ext}^{i-1}_G(U,V^{\vee\vee}/V)$.

*Proof.* Apply Theorem A and Lemma 4 to (†). □

**Corollary 3.** *Let $U$ be an arbitrary module, and $V$ a finitary module.*
*a. If $U=[U,G]$, then localization* $\operatorname{Ext}^1_G(U,V)\to\operatorname{Ext}^1_{\mathcal{L}}(U,V)$ *is injective.*
*b. If* $\operatorname{Hom}(G,k)=0$, *then localization* $\mathrm{H}^2(G,V)\to\mathrm{H}^2(\mathcal{L},V)$ *is injective.*

*Proof.* If $U=[U,G]$ then $\operatorname{Ext}^0_G(U,V^{\vee\vee}/V)=\operatorname{Hom}_{kG}(U,V^{\vee\vee}/V)=0$, whence we obtain part $a$. If $\operatorname{Hom}(G,k)=0$ then $\operatorname{Ext}^1_G(k,V^{\vee\vee}/V)=\operatorname{Hom}(G,V^{\vee\vee}/V)=0$, which yields part $b$. □

**Corollary 4.** [M, Thm 1] *If $V$ is locally completely reducible, then $[V,G]$ is completely reducible.*

*Proof.* Let $W<[V,G]$. We show that $W$ is a direct summand. Set $U=[V,G]/W$. Since $[V,G]=[V,G,G]$, $\operatorname{Hom}_G(U,W^{\vee\vee}/W)=0$. Now apply Corollary 3$a$. □

**Corollary 5.** [M, Thm 2] *If* $\mathrm{H}^1(L,V)$ *and* $\mathrm{C}_V(L)$ *are finite-dimensional for every $L$ in $\mathcal{L}$, then localization* $\mathrm{H}^2(G,V)\to\mathrm{H}^2(\mathcal{L},V)$ *is injective.*

*Proof.* If $V$ is infinite dimensional, apply Corollary 2. If not, the exact sequence

$$\cdots\to V/\mathrm{C}_V(L)\to\mathrm{H}^1(L,\mathrm{C}_V(L))\to\mathrm{H}^1(L,V)\to\cdots.$$

tells us that $\mathrm{H}^1(L,k)=0$, since otherwise $\mathrm{H}^1(L,V)$ would be infinite dimensional. Hence $\mathrm{H}^1(G,k)=0$. Now apply Corollary 3$b$. □

**Corollary 6.** [M, Thm 3] *If $V\hookrightarrow W$ is a locally split extension of $V$ by a trivial $G$-module, then there is a canonical injection* $W/\mathrm{C}_W(G)\hookrightarrow[V^\vee,G]^\vee$.

*Proof.* Note that $[V^\vee,G]^\vee=V^{\vee\vee}/\mathrm{C}_{V^{\vee\vee}}(G)$, and apply Theorem A. □

## Two Examples

For the first example, let $\Omega$ be a set of some infinite cardinal $\aleph$, let $G$ be the finitary symmetric group on $\Omega$, let $k$ be of characteristic 2, and let $V$ be the natural permutation module $k\Omega$. The Eckmann-Shapiro Lemma (see [B, Cor 2.8.4]) tells us that $\mathrm{H}^1(G,V^\vee)=\mathrm{H}^2(G,V^\vee)=k$. Now $V$ embeds in $V^\vee$, with $G$-trivial quotient (it is a quotient of $V^{\vee\vee}/V$). The associated long exact sequence contains the fragment

$$\cdots\to\mathrm{H}^1(G,V^\vee)\to\mathrm{H}^1(G,V^\vee/V)\to\mathrm{H}^2(G,V)\to\mathrm{H}^2(G,V^\vee)\to\cdots.$$

Thus $\dim_k\mathrm{H}^2(G,V)=\dim_k\mathrm{H}^1(G,V^\vee/V)=2^\aleph$. On the other hand, if we express $G$ as a union of finite symmetric groups, then we have that $\mathrm{H}^2(\mathcal{L},V)=k$. To see this, let $\Gamma$ and $\Gamma'$ be cofinite subsets, with $\Gamma\supset\Gamma'$. Set $F=\mathrm{C}_G(\Omega)$, $F'=\mathrm{C}_G(\Omega')$, and $\Delta=\Omega-\Omega'$. We have that $\mathrm{H}^2(F,V)=k\oplus\mathrm{H}^2(F,k\Delta)\oplus\mathrm{H}^2(F,k\Omega')$, while $\mathrm{H}^2(F',V)=k\oplus\mathrm{H}^2(F',k\Omega')$. Thus, taking successively smaller $\Gamma'$, we see that none of the classes in $\mathrm{H}^2(F,V)$ that take values in $k\Gamma$ survive in the limit.

Although very little of the 2-cohomology of the symmetric group is detected locally, all of the 2-cohomology for the alternating subgroup is detected locally, by Corollary 3. Moreover, all of the 2-cohomology for the symmetric group is detected locally when we replace $k$ by a field of any other characteristic.

For the second example, let $k$ have odd characteristic, let $V$ be an infinite-dimensional $k$-space with a nondegenerate symplectic form, and let $G$ be the finitary symplectic group. The form gives an embedding of $V$ in $V^\vee$. Again, the quotient $V^\vee/V$ is $G$-trivial. The extension $V \hookrightarrow V^\vee$ is nonsplit since $\mathrm{C}_{V^\vee}(G) = 0$. However, if we let $\mathcal{L}$ be the collection of symplectic groups of the finite-dimensional, nondegenerate subspaces of $V$, then we find that the extension is locally split. (See [CPS, Tbl 4.5, pp 186–187].) By Corollary 5, any $G$-trivial extension of $V^\vee$ is split.

## Skew Fields

Of the three theorems in [M], Meierfrankenfeld assumes that $k$ is commutative only in the first. The results above can be extended to the case of a skew field $k$, by carefully tracking whether $k$ should multiply from the left or right — or both. The guiding principle is the duality of Lemma 3, and so we must ensure that all the homology and cohomology groups are $k$-spaces. For Lemma 1 we need only take $X$ to be a right $k$-space and $Y$ a left $k$-space. However, for Lemma 3 and Theorem A we must further assume that $X$ is a $k$-$kG$-bimodule. This is tantamount to requiring that the $G$-action on $X$ be defined over the center of $k$. Note that Lemma 3 no longer follows from [B, Prop 2.8.5] when $k$ is not commutative. However, the direct proof given above works in general, if assume that $X$ is a bimodule.

In the applications to finitary modules, we must assume similarly that $U$ is a bimodule. In most of the corollaries this causes no difficulty. In particular, since we are taking $U = k$ in the application to [M, Thm 2] and [M, Thm 3], we find that these are subsumed by Theorem B even when $k$ is not commutative. However, in trying to apply Theorem B to [M, Thm 1], we find that we must assume that the $G$-action on $V$ is defined over the center of $k$ — which in effect echoes Meierfrankenfeld's hypothesis for this theorem.

Recently — and independently — Wehrfritz proves Corollary 4 assuming only that $k$ is finite-dimensional over its center [$W_1$]. In fact Theorem A can be extended to this situation, by tensoring with a splitting field. If $\zeta$ is the center of $k$, and $\mu$ is a maximal commutative subfield of $k$, then $\mathrm{Hom}_k(X, Y^\vee) \otimes_\zeta \mu = \mathrm{Hom}_\mu(X, Y^\vee)$. So, $\mathrm{Ext}^*_G(X, Y^\vee) \otimes_\zeta \mu = \mathrm{Ext}^*_{\mathcal{L}}(X, Y^\vee) \otimes_\zeta \mu$, whence $\mathrm{Ext}^*_G(X, Y^\vee) = \mathrm{Ext}^*_{\mathcal{L}}(X, Y^\vee)$, since field extension is faithfully flat.

The example in [$W_2$] shows that, without some hypothesis on the field, Corollary 4 would be false. (See also [ShW, Ex 1.18].)

## Acknowledgment

It is a pleasure to record my thanks to Ulrich Meierfrankenfeld, for his many insights, suggestions, and patient explanations; and also to the Math Department at Michigan State University, for their generous hospitality.